\documentclass{amsart}

\newtheorem{theorem}{Theorem}[section]

\newtheorem{corollary}[theorem]{Corollary}

\theoremstyle{definition}
\newtheorem{definition}[theorem]{Definition}

\theoremstyle{remark}

\numberwithin{equation}{section}

\newcommand{\bbC}{{\mathbb C}}
\newcommand{\bbZ}{{\mathbb Z}}
\newcommand{\bbR}{{\mathbb R}}
\newcommand{\bbP}{{\mathbb P}}
\newcommand{\calQ}{{\mathcal Q}}
\newcommand{\calZ}{{\mathcal Z}}
\newcommand{\calP}{{\mathcal P}}
\newcommand{\calH}{{\mathcal H}}
\newcommand{\calE}{{\mathcal E}}
\newcommand{\calF}{{\mathcal F}}
\newcommand{\calC}{{\mathcal C}}
\newcommand{\frakg}{{\mathfrak g}}
\DeclareMathOperator{\tr}{Tr}

\newcommand{\ol}{\overline}
\newcommand{\dbar}{\bar\partial}
\newcommand{\tD}{{\tilde\Delta}}

\begin{document}

\title{Almost Complex Structures and Geometric Quantization} 
\date{August 16, 1996}
\author{David Borthwick}
\address{Mathematics Department\\
University of Michigan\\Ann Arbor, Michigan 48109}
\curraddr{Mathematics Department\\
University of California\\Berkeley, CA 94720}
\email{borth@math.lsa.umich.edu}
\thanks{First author supported in part by NSF grant DMS-9401807}
\author{Alejandro Uribe}
\address{Mathematics Department\\
University of Michigan\\Ann Arbor, Michigan 48109}
\email{uribe@math.lsa.umich.edu}
\thanks{Second author supported in part by NSF grant DMS-9623054.}

\subjclass{Primary 53C15, 81S10}
\begin{abstract}
We study two quantization schemes for compact symplectic manifolds with almost
complex structures.  The first of these is the Spin$^c$ quantization.
We prove the analog of Kodaira vanishing for the Spin$^c$ Dirac operator,
which shows that the index space of this operator provides
an honest (not virtual) vector space semiclassically.  
We also introduce a new quantization scheme, based on a rescaled Laplacian,
for which we are able to prove strong semiclassical properties.  The two
quantizations are shown to be close semiclassically.
\end{abstract}     
\maketitle
\tableofcontents

\newcommand{\ltk}{L^{\otimes k}}
\newcommand{\eltk}{\calE\otimes\ltk}
\newcommand{\cinf}{{C^{\infty}}}
\newcommand{\ind}{\mathop{\hbox{ind}}}
\newcommand{\Hak}{\calH}
\newcommand{\Pak}{\Pi}
\newcommand{\Hs}{\calQ}
\newcommand{\Ps}{\Theta}
\newcommand{\hlap}{\Delta^\bullet}

\section{Introduction}

In the theory of geometric quantization of B. Kostant and
J-M. Souriau a key role is played by the concept of polarization.
Recall that if $X$ is a symplectic manifold, a polarization
of $X$ is a Lagrangian sub-bundle, $\calP$, of $TX\otimes\bbC$ that
satisfies the Frobenius integrability condition.
If $\calP$ is the $(0,1)$ bundle associated to an almost complex
structure, then by the Newlander-Nirenberg theorem the integrability 
condition ensures that the almost complex structure comes from
a complex structure, i.e. $X$ is a K\"ahler manifold.  The quantization of
K\"ahler manifolds has had numerous applications, from 
Bargmann's work in quantum mechanics to Kostant
and Kirillov's in the orbit method of representation theory.
In addition, this theory has excellent semi-classical properties:
see \cite{BMS}, \cite{BPU1} and \cite{BPU2}.

\medskip
Compact K\"ahler manifolds have been, until recently, the only compact
symplectic manifolds for which one has had a general method of
quantization, where by this term we mean a method for associating to
$X$ a Hilbert space $\calH$, and to functions on $X$ operators on
$\calH$.  (There exist other approaches to quantization; we mention
M.\ Karasev's \cite{Ka}.)  Recently Michele Vergne has
generalized the Kostant-Souriau scheme by replacing the polarization
with an operator of Dirac-type \cite{Ve}. 
Victor Guillemin, \cite{Gu}, has pointed out 
that a natural choice for this operator is the Spin$^c$ Dirac
operator, which is a generalization of the operator $\dbar+\dbar^*$ of the
K\"ahler theory and is constructed starting with an
almost complex structure which may not be integrable (more details
below).   The resulting ``Spin$^c$ quantization'' has attracted a great
deal of attention recently because of its nice properties with regard to
symplectic reduction \cite{DGMW}, \cite{M1}, \cite {M2}, \cite{TZ}, 
\cite{V2}.  The original definition associates a virtual vector space to $X$, 
but we will prove a vanishing theorem that shows that semiclassically this
is an ``honest'' Hilbert space.

We will also propose here another quantization method which is
well-defined semiclassically, and is based on some old results
of Victor Guillemin and one of us \cite{GU}.  Instead of generalizing
the operator $\dbar +\dbar^*$, our idea is to generalize
the Hodge Laplacian $\hlap = \dbar^*\dbar$ on sections of a holomorphic
pre-quantum
line bundle $L\to X$.  As we will see, this can also be done starting
with an almost complex structure on $X$.  
It turns out that the
semiclassical theory of this quantization method, which we'll refer to as
``almost K\"ahler quantization,'' is as good as that of K\"ahler
quantization.   This is based on the results of
\cite{GU} and the calculus of Hermite Fourier integral operators.  In fact the
methods of proof of the main results of \cite{BMS}, \cite{BPU1} and \cite{BPU2}
generalize directly to this setup.  We will also prove that semiclassically
the almost K\"ahler quantization is close to the Spin$^c$ quantization.

Before describing these two approaches in more detail, we note that 
in both schemes the quantization of functions can be defined
in the same way by the ``Toeplitz recipe.''  In either case the space of
sections has an $L^2$ inner product.   If $\Pi$ denotes the orthogonal
projection from the $L^2$ space of sections onto the quantizing Hilbert
space, $\calH$, then the quantization of $f\in\cinf(X)$ is the operator 
$\Pi M(f)$ restricted to $\calH$, where  $M(f)$ denotes the operator of
multiplication by $f$.

The Spin$^c$ Dirac operator arises as follows.
Suppose for a moment that $X$ is a compact complex manifold.  Given a
Hermitian holomorphic line bundle $L\to X$, consider the twisted Dolbeault
complex: 
$$
0\longrightarrow \cinf(X, \calE^0\otimes L)
\,\smash{\mathop{\longrightarrow}\limits^{\dbar}}\,
\cinf(\calE^1\otimes L) \longrightarrow \ldots \longrightarrow 
\cinf(\calE^n\otimes L) \longrightarrow 0,
$$
where $\calE^q$ is the bundle of forms of type $(0,q)$ and $n$ is half the
real dimension of $X$.  
The Kodaira vanishing theorem says that for $L$
sufficiently positive this complex has non-zero cohomology only in 
degree zero.  Let $\calE$ be the direct sum over $q$ of the bundles
$\calE^q$.  The vanishing theorem can also be stated in terms of the operator
$\dbar + \dbar^*$, acting on sections of $\calE\otimes L$:  If $L$ is
sufficiently positive, then the kernel of $\dbar + \dbar^*$ has only
zero-degree components (just holomorphic sections of $L$).  The index of
the Dolbeault complex is the Riemann-Roch number $RR(X,L)$, which by the
Riemann-Roch-Hirzebruch index theorem equals
$$
RR(X, L) = \int_X \hbox{Td}(X) \hbox{ch}(L)
$$
where $\hbox{Td}(X)$ is the Todd class of $X$.

Now suppose that $X$ has merely an almost complex structure 
and $L$ is a Hermitian line bundle with compatible Hermitian connection.  
The Dolbeault complex may be replaced with the ``rolled up'' version
\begin{equation}\label{rdolb}
(\dbar + \dbar^*): \cinf(X, \calE^\pm\otimes L) \to
\cinf(X, \calE^\mp\otimes L),
\end{equation}
where $\calE^+$ and $\calE^-$ are the direct sum of the $\calE^q$ for 
even and odd $q$, respectively.  The Riemann-Roch formula still computes
the index of this two-term complex.

The Spin$^c$ Dirac operator $D$ is an alternative to the operator 
$\dbar + \dbar^*$ in this situation.    For $X$ symplectic, to define it
requires choosing a compatible almost complex structure on $X$ and a
Hermitian connection on the dual of the canonical line bundle of $X$.  
(For details, see J.J Duistermaat's excellent book, \cite{D}.)
$D$ acts on the same rolled-up complex (\ref{rdolb}).
Moreover, $D$ has the same principal symbol as $2(\dbar + \dbar^*)$, 
so its index is
the again $RR(X,L)$.  But $D$ has more natural properties in many respects.  
For instance, there is an explicit local formula for the integrand in the
Riemann-Roch formula in terms of heat kernels \cite{D}.  

Following the ideas of Vergne and Guillemin, 
if $L$ is a pre-quantum line bundle of $X$ (i.e. $L$ is Hermitian with a
compatible connection and its curvature is the symplectic form) one 
defines the quantization of $X$ to be the virtual vector space 
\[
\ind D \, =\, \ker D^+ - \ker D^-,
\]
More generally, if we consider the tensor power
$\ltk$ (which has curvature $k\omega$), then $1/k$ plays the role
of Planck's constant.   One of the main 
points of this paper is to prove an analog of Kodaira vanishing for the
Spin$^c$ Dirac operator.  
Namely, we will prove in Section 2 that $\ker D^- = 0$ for
$k$ sufficiently large.   

In the almost K\"ahler quantization, which will be introduced fully in
Section 3, the quantum
Hilbert space will be defined as the span certain eigenfunctions of a rescaled
Laplacian operator.  The motivation is as follows.  In the K\"ahler setting,
the space of holomorphic sections of $\ltk$ is just the kernel of
the Hodge Laplacian $\dbar^*\dbar$.  Another natural choice of operator in this
setting would be the metric Laplacian on sections, $\Delta_k$, defined
using the Hermitian connection on $L$.  
The relationship between the two is 
\begin{equation}\label{twolap}
4 \dbar^*\dbar = \Delta_k - nk,
\end{equation}
where $2n$ is the real dimension of $X$, so the Hodge Laplacian be seen as a
rescaling of $\Delta_k$. 

In the non-integrable case, we take the right-hand side of (\ref{twolap}) as the
definition of the  rescaled Laplacian, $\hlap_k$, and we define the Hilbert
space $\Hak_k$ to be the span of the first
$d_k$ eigenfunctions of $\Delta^\bullet_k$, where 
$$
d_k\, :=\, RR(X,\ltk)\,.
$$
This is in fact a more natural generalization of the K\"ahler situation
then it might seem.   By the results of \cite{GU} and our Theorem
\ref{projdiff}, the first $d_k$ eigenvalues of
$\hlap_k$ can be bounded independently of $k$, while the rest of the spectrum 
drifts to the right with a gap of $O(k)$.  So as $k\to\infty$ the first $d_k$
eigenvalues (which are zero in the integrable case) 
remain within some fixed interval around zero even in the non-integrable case.

In fact, the methods used to prove this fact in \cite{GU} imply 
more than this.   Let $\Pak_k: L^2(X, \ltk)\to \Hak_k$ be the orthogonal
projection.  Then \cite{GU} shows that the $\Pi_k$
are the components of a projector of Szeg\"o type.  That is, they define
a generalized Toeplitz structure in the sense of \cite{BG}.  
This means that the semiclassical theorems of 
\cite{BMS}, \cite{BPU1}, and \cite{BPU2}, which are based on \cite{BG}, 
can immediately be extended to this quantization.

In Section 4 we will apply the vanishing theorem of Section 2 to 
demonstrate that the spaces $\ker D^+$ and $\Hak_k$ approach each other
semiclassically.
This gives a relation between the two methods of quantization 
(Spin$^c$ and almost K\"ahler) that is strong enough to extend some, 
but not all, of the semiclassical theorems to the Spin$^c$ quantization.

\newcommand{\inm}{\mathop{\hbox{int}}}
\newcommand{\norm}[1]{\left\Vert #1 \right\Vert}

\section{The Spin$^c$ Dirac Vanishing Theorem}
\subsection{Drift of the Laplacian}\label{drift}

Our proof of the vanishing theorem is based on a theorem 
concerning the large $k$ behavior of the spectrum of the Laplacian
acting on sections of $\eltk$, where $\calE$ is any Hermitian vector
bundle and $L$ a Hermitian line bundle, both with Hermitian connections.  In
particular, we will prove in this subsection that as $k\to\infty$ the spectrum
of $\Delta_k$ drifts to the right, at a rate governed by the
``non-flatness'' of the connection on $L$, provided the curvature of this
connection has constant rank.  For the Laplacian acting purely on sections of
$\ltk$ this fact was proven in \cite{GU}.

Start with $X$, a compact manifold of real dimension $2n$ with Riemannian
structure $\beta$.   Let $\nabla^L$ and $\nabla^\calE$ denote the 
connections on $L$ and $\calE$, respectively, and let $\nabla_k$ be the
induced connection on $\eltk$. 
The metric Laplacian on sections of $\eltk$ is $\Delta_k :=
\nabla_k^*\nabla_k$.  

Let $\omega = \hbox{curv}(\nabla^L)$, a possibly degenerate 2-form on $X$.  
For $x\in X$ define the skew-symmetric linear map $J_x: T_xX \to T_xX$ by
$$
\omega(v,w) = \beta_x(v, J_x w), \quad\hbox{for }v,w\in T_xX.
$$
The eigenvalues of $J$ are purely imaginary.  Let $\tau(x) = \tr^+ J_x :=
\mu_1+\ldots+\mu_l$, where $i\mu_j$, $j = 1,\ldots, l$, are the eigenvalues of
$J_x$ for which $\mu_j>0$.  Let $\tau_0 = \inf \tau(x)$.  

\begin{theorem}\label{driftthm}
With the definitions above, and provided the rank of $\omega$ is constant over $X$,
there exists a constant $C$ such that for all $k$ the spectrum of $\Delta_k$
lies to the right of $k\tau_0 - C$.
In particular, if $\omega, J$, and $\beta$ define compatible
symplectic, almost complex, and Riemannian structures, respectively, then
the spectrum of $\Delta_k$ lies to the right of $kn - C$.
\end{theorem}

\begin{proof}
$\calE$ will always be an induced vector bundle for some
principal $G$-bundle $\calF$, where $G$ is a compact Lie group 
(we could always take the unitary frame bundle of $\calE$, for example).  Let
$E$ be the complex vector space on which $G$ acts so that
$\calE = \calF \times_G E$.  We will also have a $\frakg$-valued connection
1-form $\vartheta$ on $\calF$ which induces $\nabla^\calE$.

A Riemannian structure can be defined on $\calF$ by first choosing an 
Ad-$G$ invariant inner product $\beta^{\frakg}$ on $\frakg$ and then defining
$$
\beta^\calF(X, Y) := \beta(d\pi(X), d\pi(Y)) + \beta^{\frakg}(\vartheta(X),
\vartheta(Y)),
$$
where $\pi:\calF\to X$ is the projection.
This $\beta^{\frakg}$ can be chosen so that the natural identification
of sections of $\eltk$ with $G$-equivariant functions $\calF \to E$ 
extends to an isomorphism
$$
L^2(X, \eltk) \cong (L^2(\calF, \pi^*\ltk)\otimes E)^G  
$$
Moreover, the action of $\Delta_k$ on the left-hand side corresponds, 
under this isomorphism, to the restriction of the operator
$$
\Delta_k^\calF\otimes I + I\otimes \mbox{Cas}_E
$$ 
to $G$-invariant sections, where $\Delta_k^\calF$ is the Laplacian on
sections of $\pi^*\ltk$ and Cas$_E
\in$ End$(E)$ denotes the Casimir operator determined by $\beta^{\frakg}$.  
The proof thereby reduces to showing that the spectrum of
$\Delta_k^\calF$ lies to the right of $k\tau_0 - C$ for $C$ independent of $k$.

Now we essentially apply the same trick to $\pi^*\ltk \to\calF$.
Let $P$ be the principal $S^1$-bundle associated with $L$, and $\alpha$ the
connection form on $P$ which induces $\nabla^L$.   And let
$\calZ = \pi^*P$, the pullback of $P$ to $\calF$:
$$
\begin{array}{ccc}
\calZ & \to &P \\
\downarrow &&  \downarrow\\
\calF &\smash{\mathop{\to}\limits^{\pi}}& X \\
\end{array}
$$
(As a bundle over $X$, $\calZ$ is just the associated principal bundle to
$\calE\otimes L$.)   We can identify sections of $\pi^*\ltk$ 
with the $k$-th isotype of the $S^1$ action on $\calZ$, 
$$
L^2(\calF, \pi^*\ltk) \cong L^2(\calZ)_k.
$$
Under this identification $\Delta^\calF$ can be written as the restriction
of an operator $\tD$ on $\calZ$ which is independent of $k$.  
(This $\tD$ is in fact the horizontal Laplacian for the bundle $\calZ \to
\calF$.) 

Let $\rho:\calZ\to X$ be the projection.
We can split the tangent space $T_z\calZ$ into $H \oplus V_{S^1} \oplus
V_{G}$, where $H$ is the horizontal lift of $T_{\rho(z)}X$, and $V_{S^1}$ 
and $V_{G}$ are tangent to the $S^1$ and $G$ actions, respectively. 
Correspondingly we have $T^*_z\calZ = H^* \oplus V^*_{S^1} \oplus
V^*_{G}$.  Using the natural isomorphisms $H^* \cong T^*_{\rho(z)}X$ and
$V^*_{G} \cong \frakg^*$, the principal symbol of $\tD$ can be written
$$
\sigma(\tD)(z,\xi) = \beta^*_{\rho(z)}(\xi_H, \xi_H) + 
\beta^{\frakg *}(\xi_{\frakg^*}, \xi_{\frakg^*}),
$$
where $\xi_H$ and $\xi_{\frakg^*}$ are the $H^*$ and $V^*_{G}$ components,
respectively, of $\xi\in T^*_z\calZ$.  The subprincipal symbol of $\tD$ 
is identically zero.

On the characteristic set $\calC = \{(z,\xi)\in T^*\calZ\backslash 0: \xi_H =
\xi_{\frakg^*} = 0\}$, $\sigma(\tD)$ vanishes to second order.  At a
point $(z,\xi)\in \calC$, 
we can therefore define the Hamilton map $F_{z,\xi}$ of
$\sigma(\tD)$, a skew-symmetric linear map in $T_{z,\xi}(T^*\calZ)$.  
The restriction that $\omega$ be of constant rank allows us to apply 
Theorem 22.3.2 of \cite{H}, a Melin-type inequality.  
In our case, with zero subprincipal symbol, the result implies that if
if we choose a first
order pseudodifferential operator $A$ on $\calZ$ such that $\sigma(A)(z,\xi) 
\le \tr^+ F_{z,\xi}$ for all ${z,\xi}\in \calC$, then there exists a 
constant $C$ so that
\begin{equation}\label{Hineq}
\langle \tD f,f \rangle \ge \langle Af, f\rangle - C \norm{f}^2
\end{equation}
for all $f\in C^\infty(\calZ)$.  Here $\tr^+$ is defined exactly as for $J_x$
above.

The map $F_{z,\xi}$ depends only on the Hessian of $\sigma(\tD)$, and one
quickly sees that to compute the non-zero eigenvalues one need only consider
the restriction of the Hessian to the symplectic subspace $H\oplus
H^* \subset T_{z,\xi}(T^*\calZ)$.  A computation in local coordinates
shows that
$$
\tr^+ F_{z,\xi} = \tau(\rho(z)) \cdot \xi_\theta,
$$
where $\xi_\theta$ is the $V^*_{S^1}$ component of $\xi$.  Let $D_\theta$ be the
generator of the $S^1$ action on $\calZ$, whose symbol is $\xi_\theta$.
We apply (\ref{Hineq}) with $A = \tau_0 D_\theta$ to give 
$$
\langle\tD f,f\rangle \ge \tau_0\langle D_\theta f, f\rangle - 
C \norm{f}^2.
$$
Noting that if $f \in C^\infty(\calZ)$ comes from a section of 
$\pi^*\ltk$ then $D_\theta f = kf$, the proof is complete.
\end{proof}

\subsection{Vanishing theorem}

Now we specialize to the case where $\beta,
\omega, J$ are compatible Riemannian, symplectic, and almost complex structures
on $X$.  As in the introduction let $\calE^q$ be the bundle of $(0,q)$-forms on
$X$, with $\calE$ denoting the direct sum over $q$ of these and $\calE^\pm$
the even and odd subbundles.  

We also choose a Hermitian connection on the dual of the canonical bundle
of $X$.  These data then define the Spin$^c$ Dirac operator $D$ \cite{D},
which we decompose into even and odd components: 
$$
D^\pm:\cinf(X, \calE^\pm\otimes \ltk) \to \cinf(X, \calE^\mp\otimes \ltk).
$$

The principal symbol of $D$ (which is also the principal of
$2(\dbar + \dbar^*)$) is given as follows.  To a vector $\xi \in T^*_xX$ we can
associate an endomorphism of $\calE_x$ by
$$
c(\xi): \nu \mapsto (\xi - iJ\xi)\wedge\nu - \inm(\beta_x^{-1}(\xi))
\nu,\qquad\hbox{for }\nu\in T^*_xX^{(0,*)},
$$
where $\inm(\cdot)$ denotes interior multiplication by a vector in $T_xX$.
This map $c$ in fact extends to an isomorphism of the complexified Clifford
algebra of $T^*_xX$ with End$(\calE_x)$.  The principal symbol of $D$ is 
given by $\sigma(D)(x,\xi) = ic(\xi)$.

Our proof of the vanishing theorem rests on the following calculation of 
$D^2$.  
\begin{theorem}\label{Dformula}
(Theorem 6.1 of \cite{D})
$$
D^2 = \Delta_k + k\sigma + R, 
$$
where $\Delta_k$ is defined as in Section \ref{drift} and $R \in
End(\calE)\otimes I$.  The operator $\sigma$ is defined by
$$
\sigma = -i \sum_{j>k} \omega(v_j, v_k) c(\xi_j) c(\xi_k),
$$
where $v_j$ is a local orthonormal frame for $TX$ and $\xi_j$ is the
corresponding dual frame.  
\end{theorem}
Note that (in contrast to $(\dbar+\dbar^*)^2$ in the K\"ahler case) $D^2$
preserves degree only mod 2 in general.  A simple calculation shows that on
forms of degree $q$, $\sigma$ acts as  multiplication by $2q-n$.  

We can now prove the main result of this section.  Our proof is inspired by
the proof in \cite{BGV} of the Kodaira vanishing theorem 
based on a Lichnerowicz-type formula.

\begin{theorem}\label{vanish}
There exist constants $C, K$ such that, for 
$\phi\in\cinf(X, \eltk)$, $k>K$, $D\phi = 0$ implies that 
$$
D\phi = 0 \quad \Longrightarrow\quad \norm{\psi} < Ck^{-1} \norm{\phi_0}, 
$$
where $\phi = \phi_0 + \psi$ is the decomposition of $\phi$ into zero and
higher degree components.
\end{theorem}
\begin{proof}
In the formula of Theorem
\ref{Dformula}, only the $R$ term will mix zero-degree components with 
those of higher degree.  We therefore have
$$
\langle \psi, D\phi\rangle = \langle \psi, (\Delta_k + k\sigma + R)\psi +
R\phi_0 \rangle.
$$
If $D\phi = 0$ then
\begin{equation}\label{psph}
|\langle \psi, (\Delta_k + k\sigma + R)\psi \rangle| = 
|\langle \psi, R\phi_0 \rangle| \le \norm{R}\>\norm{\psi}\>\norm{\phi_0}.
\end{equation}

By Theorem \ref{driftthm}, and the fact that $\sigma$ acts as
multiplication by $2q-n$, there exists a constant
$C'$ so that for large $k$
$$
|\langle \psi, (\Delta_k + k\sigma + R)\psi \rangle| >
C'k \norm{\psi}^2.
$$
Combining this with (\ref{psph}), factoring out a $\norm{\psi}$, and
absorbing $\norm{R}$ into the constant yields
$$
\norm{\psi} < Ck^{-1}\norm{\phi_0}.
$$
\end{proof}

So if $k$ is large, sections in $\ker D$ are dominated by their zero-degree components. 
If $\phi\in \ker D^-$, then $\phi_0 = 0$, so the vanishing result follows immediately:
\begin{corollary}
For $D$ acting on $\cinf(X, \eltk)$, if $k>K$ then
$$
\ker D^- = 0.
$$
\end{corollary}

\noindent
{\em Remark.}  We have stated the vanishing result in the
context of geometric quantization.  But for the proof we do not really 
need the fact that curv$(\nabla^L)$ is symplectic, or even
non-degenerate.  The more general statement is that $\ker D^- = 0$ for any
sufficiently positive Hermitian line bundle $L$ with curvature of constant
rank.

\section{Almost K\"ahler Quantization}

\subsection{Definitions}

Let $(X,\omega)$ be a compact symplectic manifold of real dimension $2n$ 
and $L\to X$ a pre-quantum
line bundle, i.e. a Hermitian line bundle with a compatible connection,
$\nabla$, whose curvature equals $\omega$.  (This means that $[\omega/2\pi]$
must be an integral cohomology class.)
Let $J$ be a compatible almost complex structure on $X$, and define
\begin{equation}\label{3a}
\forall k\in\bbZ^{+}\quad   \hlap_k\,=\, \nabla^*_k\,\nabla_k - nk\,,
\end{equation}
where
\[
\nabla_k : C^\infty(X,\ltk)\to C^\infty(X,T^*X\otimes \ltk)
\]
is the connection on $\ltk$ induced by $\nabla$. 
The adjoint of $\nabla_k$ is defined using the Riemannian structure on $X$
defined by $\omega$ and $J$, namely
\[
\beta(u,v)\,=\,\omega(u, J(v))\,.
\]
\newcommand{\psik}{\psi^{(k)}}
\newcommand{\lamk}{\lambda^{(k)}}
Let $\{\psik_j\}$ be an orthonormal basis of $L^2(X, \ltk)$ of eigenfunctions
of $\hlap_k$:
\begin{equation}\label{3b}
\hlap_k\, \psik_j\, =\, \lamk_j\,\psik_j\,,
\end{equation}
where
\[
\lamk_1\leq\lamk_2\leq\cdots\,.
\]
Recall that the Riemann-Roch polynomial of $(X,\omega)$ is by definition
\begin{equation}\label{3c}
d_k\,=\,\int_X\,e^{k\omega/2\pi}\,\hbox{Td}(T^{0,1}X)
\end{equation}
where $\hbox{Td}(T^{0,1}X)$ is the Todd class of the $(0,1)$ tangent
bundle of $X$.  The Riemann-Roch polynomial is independent of the
choice of $J$, as all such choices are homotopic.

\medskip
We are now ready to define the Hilbert space of the almost K\"ahler quantization.
\begin{definition}
For every positive integer $k$, define
\[
\Hak_k\,=\, \hbox{\rm Span}\ \{\,\psik_j\,;\,1\leq j\leq d_k\,\}\,.
\]
\end{definition}

In the holomorphic case, i.e. if $J$ is integrable and $L$ holomorphic, 
$\hlap_k$ is a multiple of $\dbar^*\dbar$ by the Bochner identity. 
Therefore, by the Riemann-Roch theorem and the Kodaira vanishing theorem, 
there exists a $K$ such that for all $k>K$, 
\begin{equation}\label{3d}
\Hak_k\,=\, \ker \hlap_k\,=\,H^0(X, \ltk)\,.
\end{equation}
In the non-integrable case there is a vestigial form of the first equality,
namely:

\begin{theorem}\label{guthm} \cite{GU}
There exist positive constants $a,b,K$ such that for all integers $k>K$
\begin{enumerate}
\item $\forall j\in\{1,\ldots , d_k\},\quad \lamk_j\in (-a, a)$.
\item $\forall j > d_k, \quad \lamk_j\geq bk$.
\end{enumerate}
\end{theorem}
Thus for large $k$ the first  $d_k$ eigenvalues of $\hlap_k$ 
(counted with multiplicities) lie in $(-a,a)$, and the rest of the
spectrum drifts to the right with a gap of $O(k)$.

\medskip
\noindent
{\em Remark 1.}  The result of \cite{GU} is more general:  One can
start with an arbitrary Riemannian metric on $(X,\omega)$ and suitably
modify the metric Laplacian, $\nabla_k^*\nabla_k$ (again by a zeroth order
perturbation), to obtain an operator for which the previous result remains true. 
Furthermore the semi-classical results that we will state below remain valid in
this setting as well.

\smallskip
\noindent
{\em Remark 2.}   Strictly speaking, the result of \cite{GU} states that for
some integer $k_0$ and some constants $a,b,K$, if $k>K$ then 
$\lamk_j\in (-a,a)$ if
$j\le d_{k+k_0}$ and $\lamk_j>bk$ if $j>d_{k+k_0}$.
However, a consequence of Theorem \ref{projdiff} is that $k_0=0$.

We define the quantization of real-valued functions on $X$ as in the
integrable case:

\begin{definition}\label{pakdef}
For every positive integer $k$, let
\[
\Pak_k\,:\, L^2(X, \ltk)\to \Hak_k
\]
denote the orthogonal projection.
For $f\in \cinf(X)$ we define
\[
T_k(f)\,=\,\Pak_k\,M(f)\,\Pak_k,
\]
where $M(f)$ denotes the operator of multiplication by $f$.
We usually regard $\Pak_k$ as an operator on $\Hak_k$.
\end{definition}

\subsection{Semi-classical results}
As we will now see, the quantization scheme defined above has excellent
semi-classical properties.  We will briefly state two theorems; the
first is a ``deformation quantization'' result, generalizing \cite{BMS}, 
and the
second is the so-called trace formula, generalizing \cite{BPU2}.  
The proofs of these generalizations are the same as those
of the cited results, because, as we will describe below, the microlocal
structure of the family of projectors $\{\Pak_k\}$ is the same in the
integrable and non-integrable cases.

\begin{theorem}\label{defthm}
For all $f,g\in C^\infty(X)$,
\begin{enumerate}
\item $\norm{T_k(f)} = \norm{f}_\infty + O(1/k)$
\smallskip
\item $\|T_k(f)\,T_k(g) - T_k(fg)\|\,=\, O(1/k)$
\smallskip
\item $\|k [ T_k(f),\,T_k(g)] - T_k(\{f,g\})\|\,=\, O(1/k)$
\end{enumerate}
where $\{\cdot,\cdot\}$ denotes the Poisson bracket on $X$.
\end{theorem}

To state the trace formula, we fix a Hamiltonian function $H\in\cinf(X)$.  Let
$E^{(k)}_i$ be the eigenvalues of $T_k(H)$.  The trace formula is an asymptotic
expansion for a weighted trace of $T_k(H)$ taken near some fixed energy $E$,
which we assume to be a regular value of $H$.  The weighting is given by a test
function $\varphi$ with compactly supported Fourier transform.  
Let $\phi$ denote the Hamiltonian flow of $H$.

\begin{theorem}\label{tracef}
If the flow $\phi$ is clean on $H^{-1}(E)$, then
we have an asymptotic expansion:
\[
\sum_{i=0}^\infty \varphi(k(E^{(k)}_i - E)) \;\sim \;
\sum_{j\in J}\;
\sum_{l=0}^\infty C_{j,l} \, e^{ik\theta_j}\, k^{(d_j-1)/2-l},
\]
where $J$ indexes the connected components of the set of pairs $(x,\tau)\in
H^{-1}(E)\times \bbR$ with $\phi_\tau(x) = x$, $d_j$ denotes the dimension of
the $j$-th component, and the angles $\theta_j$ are the holonomy angles for the
closed trajectories.  
\end{theorem}

The coefficients $C_{j,0}$ can be expressed as 
integrals over the corresponding fixed point sets.  For example, the 
$\tau = 0$ component contributes 
$$
C_{0,0} = (2\pi)^{-n}\,\hat\varphi(0)\,\hbox{vol}(H^{-1}(E)).
$$ 
For the details see \cite{BPU2}.

\medskip
There are other results whose proof is based on the 
microlocal structure.    In
\cite{BMS} several additional semiclassical theorems that follow directly from
\cite{BG} are pointed out.  For instance, a relation between quantum and
classical time evolution:
$$
\norm{e^{-iktT_k(H)} T_k(f) e^{iktT_k(H)} - T_k(f\circ\phi_t)} = O(1/k).
$$
And in addition to the trace formula, \cite{BPU2} also presents an asymptotic
expansion for a localized variant of the weighted trace, 
$$
\sum_{i=0}^\infty \varphi(k(E^{(k)}_i - E))
\Psi^{(k)}_i(x_1) \ol{\Psi^{(k)}_i(x_2)},
$$ 
where $\Psi^{(k)}_i$ is the
$i$-th eigenfunction of $T_k(H)$, in terms of classical trajectories joining
$x_1$ to $x_2$.

Finally, the results of \cite{BPU1} can be also extended to the non-integrable
case.  Here families of states in
$\Hak_k$, $k \in\bbZ_+$, are associated to certain Lagrangian submanifolds of
$X$.   These states concentrate on their associated
submanifolds as $k\to\infty$, resulting in asymptotic formulas for norms and
inner  products in terms of intersections of submanifolds.

\subsection{Microlocal structure}  

We will now briefly describe the microlocal structure that gives the
immediate extension of the proofs in 
\cite{BMS}, \cite{BPU1}, and \cite{BPU2} to our case

Each of the $\Pak_k$'s has a smooth Schwarz kernel.  Singularities arise,
and hence microlocal analysis becomes relevant, when
the $\Hak_k$'s are ``rolled up'' together, as follows.
Let $P$ be the principal $S^1$-bundle associated to $L$ as in \S2.1, and
identify
$$
L^2(X, \ltk) \cong L^2(P)_k.
$$
By considering $\Hak_k$ now as a subset of $L^2(P)$, and $\Pak_k$ as a
projector on $L^2(P)$, we may define
$$
\Hak = \bigoplus_{k=0}^\infty \Hak_k,\qquad \Pak = \bigoplus_{k=0}^\infty
\Pak_k.
$$

In the K\"ahler case, $P$ can be regarded as the boundary of 
the unit disc bundle in $L^*$, $\Hak$ is the Hardy space of $P$,
and $\Pak$ is the Szeg\"o projector for this Hardy space.
Boutet de Monvel, Sj\"ostrand, and Guillemin \cite{BG}, \cite{BS}
have shown that (in the K\"ahler case) the Szeg\"o projector is a Hermite
Fourier integral operator associated to the symplectic cone
$$
\Sigma = \{(p,r\alpha_p): p\in P, r>0\} \subset T^*P,
$$
where $\alpha$ is the connection form induced on $P$, and have worked out its
symbol.  In \cite{BG} the authors in fact develop a theory of ``generalized
Toeplitz structures,'' which are projectors expressible as Hermite FIO's
associated to a particular symplectic cone and modeled microlocally on the
Szeg\"o projector.  The semiclassical results of
\cite{BMS}, \cite{BPU1}, and \cite{BPU2} are all proven using this general
theory.  The point is that the singularities of $\Pak$ come from 
large $k$ (i.e. semiclassical) behavior  (since the $\Pak_k$ were smooth
individually).  

The basis for the proof of Theorem \ref{guthm} was the following.
\begin{theorem}\cite{GU}
The full projector $\Pak$ on $L^2(P)$ occurring in the almost K\"ahler
quantization defines a generalized Toeplitz structure associated to the
symplectic cone $\Sigma$. 
\end{theorem}

Note that in addition to the semiclassical results already stated, this
allows the immediate application of many theorems of \cite{BG} on the
structure of the full projector $\Pak$.

\newcommand{\tQ}{{\tilde Q}}
\newcommand{\tR}{{\tilde R}}

\section{Semiclassical Properties of Spin$^c$ Quantization}

Let $X, \omega, J, L$ as in \S3, and let 
$\calE$ be the bundle of type $(0,*)$ forms as
in \S2.  Choose a Hermitian connection on the dual of the canonical bundle
of $X$ so that the Spin$^c$ Dirac operator $D$ is defined.
With $K$ the constant occurring in Theorem \ref{vanish}, we make the following
definitions for the Spin$^c$ quantization.
\begin{definition}
For $k>K$ define 
$$
\Hs_k = \ker D^+ \subset L^2(X, \eltk).
$$
Let 
$$
\Ps_k: L^2(X, \eltk) \to \Hs_k
$$
denote the orthogonal projection, and for $f\in\cinf(X)$ define
$$
S_k(f) := \Ps_k M(F) \Ps_k.
$$
\end{definition}

Since $\dim \Hs_k = \ind D$ for $k>K$ by the vanishing theorem, $\dim 
\Hs_k$ is given again by the Riemann-Roch polynomial $d_k$.
The point of this section is to
compare the almost K\"ahler quantization and the Spin$^c$ quantization 
semiclassically.

In order to make such comparisons, in this section we'll 
regard the almost K\"ahler Hilbert space
$\Hak_k$ as a subspace of $L^2(X, \eltk)$ containing sections with only
zero-degree components.
Similarly, we'll extend the domain of $\Pak_k$ to $L^2(X, \eltk)$.  
\begin{theorem}\label{projdiff}
$$
\norm{\Pak_k - \Ps_k} = O(1/k).
$$
\end{theorem}

\begin{proof}
We begin by showing that $\norm{(1 - \Pak_k)\Ps_k} < Ck^{-1}$.
Given any $\phi\in L^2(X, \eltk)$, we can decompose $\phi = \eta_0 + \eta_1 +
\psi$,  where $\eta_0 \in \Hak_k$, $\eta_1$ is in the orthogonal complement of 
$\Hak_k$ in $L^2(X, \calE^0\otimes \ltk)$, and $\psi$ is the sum of
higher-degree components.  For $D\phi = 0$, $(1-\Pak_k)\Ps_k\phi = \eta_1 +
\psi$.  So we are trying to show that $D\phi = 0$
implies a bound on $\eta_1$ and $\psi$ relative to $\phi$.

Theorem \ref{vanish} already implies that for large $k$, 
$\norm{\psi} < Ck^{-1} \norm{\phi}$ for any $\phi\in\ker D$.  To show that
$\eta_1$
is small, we appeal again to Theorem \ref{Dformula}, the formula for $D^2$.  
If $D\phi = 0$ then
\begin{equation}\label{etaD}
0 = \langle \eta_1, D^2\phi\rangle = \langle \eta_1,
(\Delta_k - nk + R)(\eta_0+\eta_1) \rangle + \langle \eta_1,
R\psi \rangle
\end{equation}
Since $\eta_1$ is in the complement of $\Hak_k$, for $\Delta_k - nk$ acting
on $\eta_1$ we have a lower bound which increases with $k$:
$$
\langle \eta_1, (\Delta_k - nk + R)\eta_1) \rangle > Ck \norm{\eta_1}^2.
$$
On the other hand, the action of $\Delta_k - nk$ on $\Hak_k$ is bounded
uniformly in $k$, so
$$
\langle \eta_1, (\Delta_k - nk + R)\eta_0) \rangle < C \norm{\eta_0}
\norm{\eta_1}.
$$
Returning to (\ref{etaD}) and using the bound on $\psi$, 
we obtain
$$
k \norm{\eta_1}^2 <  C \norm{\eta_0} \norm{\eta_1} + Ck^{-1} \norm{\eta_1}
\norm{\phi}, 
$$
which leads directly to the desired $\norm{\eta} < Ck^{-1} \norm{\phi}$.

This completes the proof that $\norm{(1 - \Pak_k)\Ps_k} < Ck^{-1}$.  Noting
that $\Pak_k(1 - \Ps_k) = ((1 - \Ps_k)\Pak_k)^*$, we can use
the same argument to show that $\norm{\Pak_k(1 - \Ps_k)} < Ck^{-1}$.
Since $\Pak_k(1 - \Ps_k) - (1 - \Pak_k)\Ps_k = \Pak_k - \Ps_k$, this completes
the proof. 
\end{proof}

We return now to clarify the remark made after Theorem \ref{guthm}.  
The version of this theorem proven in \cite{GU} involves
$d_{k+k_0}$ eigenvalues, where $k_0$ is some fixed but unknown integer,
instead of $d_k$.
So we should really have defined $\dim \Hak_k$ as the span
of the first $d_{k+k_0}$ eigenfunctions.
But a simple consequence of Theorem \ref{projdiff} is that with $\Hak_k$ so
defined, $\dim \Hak_k = \dim \Hs_k$ for $k$ sufficiently large (large
enough so that $\norm{\Pak_k - \Ps_k} < 1$).  
So in fact $k_0 = 0$.

\bigskip
Theorem \ref{projdiff} has some direct applications to semiclassical
analysis of  ``Spin$^c$ Toeplitz operators.''  

\begin{theorem}\label{asymthm}
For $f\in \cinf(X)$,
$$
\norm{S_k(f) - T_k(f)} = O(1/k).
$$
\end{theorem}

\begin{corollary}\label{asymcor}
For all $f,g \in \cinf(X)$,
\begin{enumerate}
\item $\norm{S_k(f)} = \norm{f}_\infty + O(1/k)$
\smallskip
\item $\norm{S_k(f)S_k(g) - S_k(fg)} = O(1/k)$
\end{enumerate}
\end{corollary}

Theorem \ref{asymthm} follows immediately from Theorem \ref{projdiff}, and
the corollary follows because of Theorem \ref{defthm}.

\bigskip
Some additional results may be derived for a more general notion of Toeplitz
operator.  
We may regard the space $\Hs_k$ as a subspace of $L^2(\calZ)$,
where $\calZ$ is the principal bundle introduced in Section 2.  So given a
pseudodifferential operator $Q$ on $\calZ$, we can define
$$
S_k(Q) := \Ps_k Q\Ps_k.
$$  

\begin{theorem}\label{comnorm}
Let $Q$ and $R$ be pseudodifferential operators on $\calZ$ of order $q$ and 
$r$, respectively.  Then
\begin{equation}\label{normq}
\norm{S_k(Q)} = O(k^{q})
\end{equation}
and
\begin{equation}\label{normc}
\norm{[S_k(Q), S_k(R)]} = O(k^{q+r-1})
\end{equation}
\end{theorem}
\begin{proof}
Let $D_\theta$ be the generator of the $S^1$
action on $\calZ$, as in Section 2.  
To prove (\ref{normq}) for $q\le 0$, note that
$(D_\theta)^{-q} Q$ is an operator of order zero, and thus bounded because
$\calZ$ is compact.  Let $\Ps = \oplus_k \Ps_k$, as a projector on
$L^2(\calZ)$.  Since $\norm{\Ps} = 1$, we have 
$\norm{\Ps (D_\theta)^{-q} Q\Ps}<C$.  So $\norm{S_k((D_\theta)^{-q} Q)}<C$
for each $k$, where $C$ is independent of $k$.  The proof is concluded by
noting that $\Ps_k D_\theta = k\Ps_k$.   For $q>0$ we must choose a
parametrix $B$ for $D_\theta$, and consider $S_k(B^qQ)$.  Then $\Ps_k B =
k^{-1} \Ps_k$ plus $k^{-1}$ times a smoothing operator, and we may use the
result for $q>0$ to control this correction.

To prove (\ref{normc}) we need to appeal to one of the basic facts about
generalized Toeplitz structures proven in \cite{BG}, Proposition 2.13.  Namely
given a pseudodifferential operator $Q$, there exists another
pseudodifferential operator $\tQ$ such that $\Pak Q\Pak = \Pak \tQ\Pak$ and
$[\Pak, \tQ] = 0$.  $\tQ$ necessarily has the same principal symbol as $Q$.

It turns out that $\tQ$ also commutes with $\Ps_k$ to leading order in $k$.  
To see this, note that $[\Ps_k,
\tQ] = [\Ps_k - \Pak_k, \tQ]$. Inserting $D_\theta^{-q}$ as above, and using
Theorem \ref{projdiff}, we obtain
\begin{equation}\label{normrq}
\norm{[\Ps_k, \tQ]} < Ck^{q-1}.
\end{equation}

Returning to the second statement of the theorem,
$$
[S_k(Q), S_k(R)] = [S_k(\tQ), S_k(\tR)] + [S_k(Q-\tQ), S_k(\tR)] + 
[S_k(Q), S_k(R-\tR)]
$$
The norms of the last two commutators on the right can be bounded by
$Ck^{q+r-1}$, by (\ref{normq}), so it only remains to bound $[S_k(\tQ),
S_k(\tR)]$.   Observe that
$$
[S_k(\tQ), S_k(\tR)] = S_k([\tQ, \tR]) + \Ps_k \tQ [\Ps_k, \tR]\Ps_k - 
\Ps_k \tR [\Ps_k, \tQ]\Ps_k.
$$
All terms on the right have norm bounded by $Ck^{q+r-1}$; the first by
(\ref{normq}) directly, the second and third by (\ref{normq}) and 
(\ref{normrq}).  
\end{proof}

Of course, by the same argument (\ref{normq}) holds for $T_k(Q)$.  And in
fact there is a refined version of (\ref{normrq}) for almost K\"ahler
Toeplitz operators.
The theory of \cite{BG} allows one to pick off the leading
term in the commutator: 
$$
\norm{[T_k(Q), T_k(R)] - k^{-1} T_k(P)} = O(k^{q+r-2}),
$$
where $\sigma(P) = \{\sigma(Q), \sigma(R)\}$.  Part 3 of Theorem
\ref{defthm} follows immediately from this fact \cite{BMS}.  Our comparison
between $\Pak_k$ and $\Ps_k$ is not strong enough to extend this result 
to the Spin$^c$ case.

We note also that Theorem \ref{asymthm} implies that 
dist(Spec$(S_k(f)),$ Spec$(T_k(f)))$ is $O(1/k)$.
This is enough to extend the Szeg\"o theorem of \S13 of
\cite{BG} to Spin$^c$ Toeplitz operators, but not enough to extend the trace
formula (Theorem \ref{tracef}).  

\section{Concluding Remarks}

We have seen that there are two methods (at least) of quantizing triples
$(X,\omega,J)$ consisting of a compact symplectic manifold with a 
compatible almost complex structure.  Although we have demonstrated that
the almost K\"ahler quantization
has very good semiclassical properties, we have yet to prove the same for 
the Spin$^c$ quantization.  This seems to require dealing with Fourier integral
operators of Hermite type acting on sections of vector bundles.  Although there
do not seem to be any serious obstacles, as always when dealing with systems
matters are more involved.  We hope to provide the details of this approach
in the future.

\medskip
From a purely geometrical point of view, the quantizations discussed in this
paper raise a number of natural questions.  Recall a classical
construction in algebraic geometry: if $L\to X$ is a holomorphic line bundle
such that the corresponding linear system is base-point free, then $L$ defines
a map of $X$ to a projective space,
\[
F\,:\,X\longrightarrow \bbP\left(H^0(X,L)^*\right)\,.
\]
The definition of $F$ is the following:
$\forall x\in X$, $F(x)$ is the hyperplane in $H^0(X,L)$ consisting of all
holomorphic sections vanishing at $x$.  (The condition of being base-point 
free means precisely that $\forall x\in X$ such a set is indeed a hyperplane.)
In the setting considered in this article, the Hilbert spaces of
both quantization
schemes are subspaces of the space of sections of a vector bundle, so one
can attempt to define a map $F$ precisely as above, with $\Hak_k$ (or
$\Hs_k$) replacing $H^0(X,L)$.
We claim that if $k$ is sufficiently large then this definition is
possible, at least for the almost K\"ahler case.
That is, for a given $X$ for all $k$ sufficiently large and for all $x\in X$
there is a
$\psi\in\Hak_k$ such that $\psi(x)\not= 0$, and therefore the space of elements
in $\Hak_k$ vanishing at $x$ has codimension one.  Questions on the geometry
of these maps and possible relations with the  Gromov-Seiberg/Witten
invariants of $X$ will be investigated in future work.


\begin{thebibliography}{99}

\bibitem{BMS}  M.\ Bordemann, E.\ Meinrenken, and M.\ Schlichenmaier:
Toeplitz quantization of K\"ahler manifolds and gl$(N)$, $N\to\infty$
limits,  {\it Comm. Math. Phys.} {\bf165} (1994) 281--296.

\bibitem{BPU1} D.\ Borthwick, T.\ Paul, and A.\ Uribe: Legendrian distributions
with applications to relative Poincar\'e series, {\it Inv. Math.} {\bf122}
(1995) 359--402.

\bibitem{BPU2} D.\ Borthwick, T.\ Paul, and A.\ Uribe: Semiclassical spectral
estimates for Toeplitz operators, preprint (1995).

\bibitem{BG} L.\ Boutet de Monvel and V. Guillemin:
{\it The spectral theory of Toeplitz operators.} 
Annals of Mathematics Studies No.\ 99, Princeton University Press,
Princeton, New Jersey (1981).

\bibitem{BGV} N. Berline, E. Getzler, and M. Vergne: {\it 
Heat kernels and Dirac operators}, Srpinger-Verlag, Berlin (1992).

\bibitem{BS} L. Boutet do Monvel and J. Sj\"ostrand, Sur la singularit\'e des
noyaux de Bergmann et de Szeg\"o, {\it Asterisque} {\bf 34--35} (1976)
123--164. 

\bibitem{D}  J.\ J.\ Duistermaat: {\it The Heat Kernel Lefschetz
Fixed Point Formula for the Spin-c Dirac Operator}, Birkh\"auser, Boston
(1996). 

\bibitem{DGMW}  H. Duistermaat, V. Guillemin, E. Meinrenken, and S. Wu:
Symplectic reduction and Riemann-Roch for circle actions, {\it
Math. Res. Lett.} {\bf 2} (1995) 259--266. 

\bibitem{Gu} V.\ Guillemin: Reduced phase spaces and Riemann-Roch, in: 
Lie theory and geometry in honour of B. Kostant, {\it Progr. Math.} {\bf123}, 
Birkh\"auser, Boston (1994) 305--334.

\bibitem{GU} V.\ Guillemin and A.\ Uribe:  The Laplace operator on the $n$-th
tensor power of a line bundle: eigenvalues which are uniformly bounded in $n$,
{\it Asymptotic Anal.} {\bf1} (1988) 105--113.

\bibitem{Ka} M. V. Karasev:  Simple quantization formula, in: Symplectic
geometry and mathematical physics (Aix-en-Provence, 1990), {\it
Progr. Math.}  {\bf99}, Birkh\"auser, Boston (1991) 234--244.

\bibitem{H} L. H\"ormander: {\it The Analysis of Linear Partial Differential
Operators III}, Springer-Verlag, Berlin (1983).

\bibitem{M1} E. Meinrenken: On Riemann-Roch formulas for multiplicities,
{\it J. Amer. Math. Soc.} {\bf9} (1996) 373--389.

\bibitem{M2} E. Meinrenken: Symplectic surgery and the Spin$^c$ Dirac
operator, {\it Adv. in Math.}, to appear (1996).

\bibitem{Sj} R.\ Sjamaar: Symplectic reduction and Riemann-Roch formulas 
for multiplicities, {\it Bull. AMS} {\bf33} (1996) 327--338.

\bibitem{TZ}  Y. Tian and W. Zhang:  Symplectic reduction and analytic
localization, preprint (1996).

\bibitem{Ve} M.\ Vergne: Geometric quantization and equivariant cohomology, 
in: First European Congress of Mathematics, vol. 1, (Paris, 1992) 
{\it Progr. Math.} {\bf 119}, Birkh\"auser, Boston (1994) 249--295.

\bibitem{V2} M. Vergne: Multiplicity formula for geometric quantization,
Parts I and II, {\it Duke. Math. J.} {\bf82} (1996) 143--179, 181--194.

\end{thebibliography}
\end{document}